\documentclass[12pt]{amsart}
\usepackage{a4wide,graphicx,amscd,amssymb}

\DeclareMathOperator\Ob{Ob}

\def\id{\mathrm{id}}

\newtheorem{theorem}{Theorem}[section]

\newtheorem{lemma}[theorem]{Lemma}
\newtheorem{proposition}[theorem]{Proposition}
\theoremstyle{definition}
\newtheorem{definition}[theorem]{Definition}
\theoremstyle{remark}
\newtheorem{remark}[theorem]{Remark}
\numberwithin{equation}{section}

\begin{document}

\title{Algebraic $K$-theory of Fredholm modules and KK-theory}
\author{Tamaz Kandelaki}
\address{Author's address: A. Razmadze Mathematical Institute, 1, M. Aleksidze St., 380093, Tbilisi, Georgia}
\email{nichbisi@yahoo.com}
\address{National Center for Science and Technology, Tbilisi, Georgia}

\subjclass[2000]{19K35, 19K33, 19D23, 19D99, 46L99}

\keywords{Nonconnective algebraic $K$-theory spectrum, category of Fredholm modules, $KK$-theory}%

\dedicatory{Dedicated to Hvedri Inassaridze}

\begin{abstract}
Kasparov $KK$-groups $KK(A,B)$ are represented as homotopy groups
of the Pedersen-Weibel nonconnective algebraic $K$-theory spectrum
of the additive category of Fredholm $(A,B)$-bimodules for $A$ and
$B$, respectively, a separable and $\sigma$-unital trivially
graded real or complex $C^*$-algebra acted upon by a fixed compact
metrizable group.
\end{abstract}

\maketitle

\section*{Introduction}

In noncommutative topology and differential geometry some of the
useful and powerful tools are methods of algebraic $K$-theory,
Kasparov's $KK$-theory, spectra and so on. Therefore a comprehensive
study of relationships between them may be considered as an
interesting task. The main goal of this paper is characterization of
Kasparov $KK$-groups as algebraic K-groups of an additive category.
On first view, calculation of the $KK$-theory by algebraic
$K$-theory seems to be highly improbable, as algebraic $K$-theory
and $KK$-theory are independent, highly nontrivial, theories, having
almost no connections with each other. The key is thus to find
suitable objects which make sense for both algebraic $K$-theory and
$KK$-theory. In this paper we concentrate on the additive
$C^*$-category $\mathrm{Rep}(A,B)$, namely the category of Fredholm
modules, where $A$ is a separable and $B$ is a $\sigma$-unital real
or complex $C^*$-algebra with action of a fixed compact second
countable group. Our main result claims the natural isomorphisms
\begin{equation*}
\mathbb{K}^a_n(\mathrm{Rep}(A;B))\simeq KK_{n-1}(A;B),
\end{equation*}
where $\mathbb{K}^a_n $ denote the algebraic $K$-functors
isomorphic to Quillen's $K$-functors in nonnegative dimensions,
and isomorphic to Pedersen-Weibel $K$-functors in negative
dimensions.

There are already several papers dedicated to interpretations of
$KK$-theory, each with their own advantages. Let us point out in
brief some fundamental papers of this sort which are sources of
further research. These are G.~Kasparov's interpretation of
$KK$-theory in terms of extensions of $C^*$-algebras \cite{kas2};
J.~Cuntz's result in \cite{cun}, based on the homotopy category of
$C^*$-algebras; N.~Higson's approach, considering $KK$-theory as
an universal enveloping additive category of the category of
separable $C^*$-algebras \cite{cun}; the interpretation of the
$KK$-theory with the aid of the category of $C^*$-algebras and
asymptotic homomorphisms, due to A.~Connes and N.~Higson
\cite{cohi}, and R. Meyer's and R. Nest's joint paper \cite{mene},
where Kasparov category $KK$ turned into a triangulated category.

Let us shortly review some known results in $KK$-theory established
by methods of topological $K$-theory and spectra.

In \cite{pwe1} E.~K.~Pedersen and C.~Weibel showed that values on finite
$CW$-complexes $X$ of the homology theory associated with the nonconnective
algebraic $K$-theory spectrum of a unital ring $R$ may be interpreted as the
algebraic $K$-groups (up to a shift in dimension) of a suitably constructed
additive category $\mathcal{C}_{O(X)}(R)$. According to this result, in
\cite{rose} J.~Rosenberg showed that
$$
K_1(\mathcal{C}_{O(X)}(R))\simeq KK(C(X),R)
$$
for a unital $C^*$-algebra $R$. He also constructed algebraic
$KK$-theory spectra, denoted by $\mathbf{KK}(A,B)$, having the
property
$$
\pi_0(\mathbf{KK}(A,B))\simeq KK(A,B).
$$
The similar question for nonzero dimensions has been left open in
that paper.

In \cite{hpr}, there have been constructed
$\mathbb{K}\mathbb{K}^{\mathrm{top}}(C_0(X),B)$-spectra, related to
Kasparov groups $KK_n(C_0(X),B)$, which were used to construct the
splitting assembly map
$$
A:\mathbb{KK}^{\mathrm{top}}(C_0(X),B)\rightarrow
\mathbb{K}^{top}(C^*(X)),
$$
where $C^*(X)$ is the $C^*$-algebra of the coarse space $X$. A
large number of results on the Novikov conjecture can be included
under this scheme \cite{hpr}.

More general approaches to non-equivariant $KK$-theory as homotopy
groups of a spectrum can be found in the following papers. These
are, \cite{hothoms}, where T.~G.~Houghton-Larsen and K.~Thomsen,
utilizing spaces of $C^*$-extensions, have constructed $KK$-theory
spectra; and \cite{mitch}, where P.~Mitchener, using methods of
topological $K$-theory, symmetric spectra and $C^*$-categories has
defined $KK$-theory spectra, too.

Let us make few remarks concerning our approach. In \cite{kan1},
\cite{kan4} we have calculated topological and Karoubi-Villamayor
$K$-groups of the $C^*$-category $\mathrm{Rep}(A,B)$  which are
related to $KK$-groups by the equalities
\begin{equation}\label{krv}
K_n^t(\mathrm{Rep}(A,B))=K_n^{\mathrm{KV}}(\mathrm{Rep}(A,B))=KK_{n-1}(A,B),
\;\;\;n\geq 0,
\end{equation}
where $A$ and $B$ are $G$-$C^*$-algebras, separable and $\sigma$-unital
respectively. Thus the additive $C^*$-category $\mathrm{Rep}(A,B))$ is a good
object for our purposes at first sight.

In the article \cite{kan2} it was announced that an isomorphism
similar to \ref{krv} for algebraic $K$-groups also holds, as well.
The present paper is an attempt to explain in detail the results
announced in \cite{kan2}. As a consequence construction of a
nonconnective algebraic $KK$-theory spectrum
$\mathbb{KK}^{\mathrm{alg}}(A,B)=\mathbb{K}(\mathrm{Rep}(A,B))$
arises, where the right hand side is the Pedersen-Weibel
nonconnective algebraic $K$-theory spectrum of the
idempotent-complete additive category $\mathrm{Rep}(A,B)$. This
spectrum has the property
$$
\pi_n(\mathbb{KK}^{\mathrm{alg}}(A,B))\simeq KK_{n-1}(A,B),
$$
for all $n\in\mathbb Z$. In a further paper we hope to apply the
algebraic $KK$-theory spectra, in particular, to problems related
to the Novikov conjecture. Our approach is mainly based on the
author's unpublished preprint \cite{kan5}.

\section{The main Theorem and an outline of the proof\label{exe}}

The purpose of this section is to summarize some of the concepts
that are needed for the formulation of the main theorem, and then to
give an outline of the proof. First, we recall the definitions of a
$C^*$-category and an idempotent-compete $C^*$-category; and in
Subsection 1.1 an idempotent-complete $C^*$-category
$\mathrm{Rep}(A,B)$ is constructed that we will need later on.

Let $A$ be a  category such that for any pair $(a,b)$ of objects in
$A$, the set $\hom(a,b)$ is equipped with the structure of a Banach
space in such a way that composition is a continuous $k$-bilinear
map. Such a category is said to be a Banach category over $k$, or
simply a Banach category. A Banach category $A$ is called a
$C^*$-category if it is equipped with a family of anti-linear maps
$*:\hom(a,b)\rightarrow\hom(b,a)$ for any $a,b\in\Ob(A)$ such that
\begin{enumerate}
\item $(f^*)^*=f$;
\item $(fg)^*=g^*f^*$, if $fg$ exists;
\item $\|f^*\|=\|f\|$;
\item $\|f^*f\|=\|f\|^2$, if $k$ is the complex numbers; and
$\|f\|^2\leq\|f^{*}f+g^{*}g\|^2$, if $k$ is the real numbers.
\item For any morphism $f:a\rightarrow b$ in $A$ the morphism $f^*f$
is a positive element of the $C^*$-algebra $\hom(a,a)$.
\end{enumerate}

Let $A$ and $B$ be $C^*$-categories. A functor $\mathcal{F}
:A\rightarrow B$ is said to be a $*$-functor if
\begin{itemize}
\item $\mathcal{F}(f+g)=\mathcal{F}(f)+\mathcal{F}(g);$
\item $\mathcal{F}(\lambda f)=\lambda \mathcal{F}(f);$
\item $\mathcal{F}(f^*)=\mathcal{F}^*$,
\end{itemize}
where $\lambda\in k$, and $f$ and $g$ are morphisms in $A$. (cf.
\cite{clr}, \cite{kan1}, \cite{kan3}).

We say that a $*$-functor is faithful if it is injective on both
objects and morphisms. Any $*$-functor is norm-nonincreasing.
Moreover, a faithful $*$-functor preserves norms \cite{clr}.

The category $\mathcal{H}(k)$ of separable Hilbert spaces and
bounded linear maps has a natural structure of a $C^*$-category.
There exists a faithful $*$-functor from every $C^*$-category into
$\mathcal{H}(k)$.

Let $A$ be a $C^*$-category and $I\subset\hom A$. Put $\hom_I(a,b)=
\hom(a,b)\cap I$. Then $I$ is called a left ideal if $\hom_I(a,b)$
is a linear subspace of $\hom(a,b)$ and $f\in\hom_I(a,b)$, $g\in
\hom(b,c)$ imply $gf\in\hom_I(a,c)$. A right ideal is defined
similarly. $I$ is a two-sided ideal if it is both a left and a right
ideal. An ideal $I$ is closed if $\hom(a,b)_I$ is closed in
$\hom(a,b)$ for each pair of objects. A closed two-sided ideal $I$
is called a $C^*$-ideal if $I = I^*$. Every $C^*$-ideal determines
an equivalence relation on the morphisms of $A$: $f\sim g$ if
$f-g\in I$; the set of equivalence classes $A/I$ can be made into a
$C^*$-category in a unique way by requiring that the canonical map
$f\mapsto\hat{f}$ gives rise to a $*$-functor $A\rightarrow A/I$.
Arguing as for $C^*$-algebras, one can show that every closed ideal
is a $C^*$-ideal \cite{clr}. Among the $C^*$-categories containing
$I$ as a $C^*$-ideal there exists a universal one, the so called
multiplier $C^*$-category. We will denote by $M(I)$ the multiplier
$C^*$-category of a $C^*$-ideal $I$ \cite{clr}.

Let $A$ be an additive category. \textit{Idempotent completion} of
$A$ is an additive category $\hat{A}$ whose objects have the form
$(a,q)$, where $a$ is an object in $A$ and $q$ is an idempotent in
$\hom(a,a)$, and a morphism $f:(a,q)\rightarrow(a',q')$ is a
morphism $f:a\rightarrow a'$ such that $fq=q'f=f$. There is a
natural functor $A\rightarrow\hat{A}$, defined by assignments
$a\mapsto(a,1_a)$ and $f\mapsto f$. An additive category $B$ is said
to be \textit{idempotent-complete} if there are an additive category
$A$, and an additive functor $F:B\rightarrow\hat{A}$ which is an
equivalence of categories.

Note that for an additive $C^*$-category $A$ the category $\hat{A}$
is not necessarily a $C^*$-category. Below we will adapt the above
construction to the case of additive $C^*$-categories.

Recall that a projection $p$ in a $C^*$-category is a morphism with
the properties $p^* = p$ and $p^2 = p$, i.~e., a projection is a
selfadjoint idempotent.

Let $A$ be an additive $C^*$-category. Consider the additive $C^*$-category
$\tilde{A}$ with objects of the form $(a,p)$, where $a\in\Ob(A)$ and
$p\in\hom(a,a)$ is a projection. A morphism from $(a,p)$ to $(b,q)$ is a
morphism $f:a\rightarrow b$ in $A$ such that $fp = qf = f$. Composition of
morphisms is the same as in $A$. The sum is given by $(a,p)\oplus(b,q) =
(a\oplus b,p\oplus q)$, and the norm of morphisms is inherited from $A$
\cite{kan1}. There is a natural functor $\nu :\tilde{A}\rightarrow\hat{A}$
defined by identity maps on objects and morphisms.

Let us show the following simple lemma.

\begin{lemma}
\label{lemaux}
Let $A$ be an additive $C^*$-category. Then $\tilde{A}$ is an
idempotent-complete $C^*$-category.
\end{lemma}

\begin{proof}
Consider the natural additive functor $\nu :\tilde{A}\rightarrow\hat{A}$,
which is, of course, faithful. Let us show that $\nu$ is a full functor.
Indeed, if $q\in\hom(a,a)$ is an idempotent then
$$
p=((2q^*-1)(2q-1)+1)^{\frac12}\;\cdot q\cdot\;((2q^*-1)(2q-1)+ 1)^{-\frac12}
$$
is a projection and the pairs $(a;q)$ and $(a;p)$ are isomorphic in $\hat{A}$
via the morphism
$$
p((2q^*-1)(2q-1)+1)^{\frac12}q.
$$
\end{proof}

Now, we define examples of additive $C^*$-categories which are
used in the remaining part of paper.

\subsection{On additive $C^*$-categories $Rep(A,B)$ and $\mathrm{Rep}(A,B)$}

Let $\mathcal{H}_G(B)$ be the additive $C^*$-category of countably
generated right Hilbert $B$-modules equipped with a $B$-linear,
norm-continuous $G$-action over a fixed compact second countable
group $G$ \cite{kas1}. Note that the compact group acts on the
morphisms by the following rule: for $f:E\rightarrow E'$ the
morphism $gf:E\rightarrow E'$ is defined by the formula
$(gf)(x)=g(f(g^{-1}(x)))$.

The category $\mathcal{H}_G(B)$ contains the class of compact
$B$-homomorphisms \cite{kas1}. Denote it by $\mathcal{K}_{G}(B)$. Known
properties of compact $B$-homomorphisms imply that $\mathcal{K}_{G}(B)$ is a
$C^*$-ideal \cite{clr} in $\mathcal{H}_G(B)$.

Objects of the category $Rep(A,B)$ are pairs of the form $(E,\varphi)$, where
$E$ is an object in $\mathcal{H}_G(B)$ and
$\varphi:A\rightarrow\mathcal{L}(E)$ is an equivariant $*$-homomorphism. A
morphism $f:(E,\phi)\rightarrow(E',\phi')$ is a $G$-invariant morphism
$f:E\rightarrow E'$ in $\mathcal{H}_G(B)$ such that
$$
f\phi(a)-\phi'(a)f\in\mathcal{K}_G(E,E')
$$
for all $a\in A$. The structure of a $C^*$-category is inherited from
$\mathcal{H}_G(B)$. It is easy to see that $Rep(A,B)$ is an additive
$C^{*}$-category, not idempotent-complete.

Now, we are ready to construct our main $C^*$-category, that is
$\mathrm{Rep}(A,B)$. Objects of it are triples $(E,\phi,p)$, where $(E,\phi)$
is an object and $p:(E,\phi)\rightarrow (E,\phi)$ is a morphism in $Rep(A,B)$
such that $p^*=p$ and $p^2=p$. A morphism
$f:(E,\phi,p)\rightarrow(E',\phi',p')$ is a morphism
$f:(E,\phi)\rightarrow(E',\phi')$ in $Rep(A,B)$ such that $fp=p'f=f$. In
detail, $f$ must satisfy
\begin{equation}
f\phi(a)-\phi'(a)f\in \mathcal{K}(E,F) \textrm{ and } fp=p'f=f.
\end{equation}
So, by definition
$$
\mathrm{Rep}(A,B)=\widetilde{Rep(A,B)}.
$$
The structure of a $C^*$-category on $\mathrm{Rep}(A,B)$ comes from the
corresponding structure on $Rep(A,B)$.

Let $S_G$ denote the category of trivially graded separable $C^*$-algebras
over $k$ with an action of the compact second countable group $G$ and
equivariant $*$-homomorphisms. Functors $\mathbb{K}_{n}^a$ are defined by
$$
\mathbb{K}^a_n(A)=\pi _n\mathbb{K}(A), \;\;\; n\in \mathbb{Z},
$$
where $\mathbb{K}(A)$ is the Pedersen-Weibel nonconnective algebraic
$K$-theory spectrum \cite{pwe} of an idempotent-complete additive category
$A$. Functors $\mathbb{K}_{n}^t$ are the topological $K$-functors on
idempotent-complete additive $C^*$-categories, defined by Karoubi
\cite{kar},\cite{kar1}. For simplicity, Kasparov's groups $KK^{-n}_G(A,B)$
will be denoted by $KK_n(A,B)$.

Now, we present our main result in the following theorem.
\begin{theorem}
\label{corollary} Let $B$ be a $\sigma$-unital trivially graded $C^*$-algebra
with an action of a second countable compact group $G$. There are natural
isomorphisms
\begin{equation}
\label{fmls} \mathbb{K}^a_n(\mathrm{Rep}(-;B))\simeq
\mathbb{K}^t_n(\mathrm{Rep}(-;B))\simeq KK_{n-1}(-;B)
\end{equation}
of functors on the category $S_G$, for all $n\in\mathbb{Z}$.
\end{theorem}

\textit{Outline of proof}. Theorem \ref{corollary} is a consequence of the
argument presented below.

A family $\mathbb{H}=\{H_n\}_{n\in \mathbb{Z}}$ of contravariant functors
from $S_G$ to the category of abelian groups and homomorphisms is said to be
a \textit{stable cohomology theory} on the category $S_G$ if
\begin{enumerate}
\item $\mathbb{H}$ \textit{has the weak excision property}. Namely,
for any exact \textit{proper} sequence
$$
0\rightarrow I\rightarrow B\rightarrow A\rightarrow 0
$$
(which means that the involved epimorphism admits an equivariant completely
positive contractive section) of algebras in $S_G$ there exists a natural
homomorphism $\delta_n:H_n(I)\rightarrow H_{n-1}(A)$, for any $n\in\mathbb
Z$, such that the resulting natural sequence of abelian groups (extending in
both directions)
$$
\cdots \rightarrow H_n(A)\rightarrow H_n(B)\rightarrow
H_n(I)\stackrel{\delta_n}{\rightarrow }H_{n-1}(A)\rightarrow
\cdots
$$
is exact.
\item $\mathbb{H}$ \textit{is stable}. This means that if
$e_A:A\rightarrow A\otimes\mathcal{K}$ is a homomorphism defined by the map
$a\mapsto a\otimes p$, where $p$ is a rank one projection in $\mathcal{K}$,
then $H_n(e_A):H_n(A\otimes \mathcal{K})\rightarrow H_n(A)$ is an
isomorphism, where $\mathcal{K}$ is the $C^*$-algebra of compact operators on
a separable Hilbert space $\mathcal{H}$ over $k$, with the trivial action of
the group $G$.
\end{enumerate}

Denote by $C_k(S^1)$ the $C^*$-algebra of continuous complex functions on the
standard unit circle $S^1$ of modulus one complex numbers, in case $k$ is the
field of complex numbers; while if $k$ is the field of reals, let it be the
subalgebra of the former consisting of functions invariant under the
conjugation defined by the map $f(z)\mapsto\overline{f(\bar{z})}$. It is
clear that any continuous complex function $f:S^1\rightarrow \mathbb{C}$ may
be represented in the form
$$
f(z)=\frac{f(z)+\overline{f(\bar{z}})}{2}+i \frac{(-i
f(z))+(\overline{-if(\bar{z})})}{2},
$$
which means that $C_{\mathbb{C}}(S^1)$ is the complexification of
the real $C^*$-algebra $C_{\mathbb R}(S^1)$.

Let $\mathbf{T}_{k}$ be the Toeplitz $C^*$-algebra---the universal
$C^*$-algebra over $k$ generated by an isometry $v$. There is a conjugation
on ${\mathbf T}_{\mathbb{C}}$ defined by the equality $\bar{v}=v$ on the
generator $v$ of ${\mathbf T}_{\mathbb{C}}$. According to the universal
property of the Toeplitz algebra, one gets a natural homomorphism ${\mathbf
T}_{\mathbb{R}}\rightarrow {\mathbf T}_{\mathbb{C}}$ so that the induced
homomorphism ${\mathbf T}_{\mathbb{R}}\;\otimes
_{\mathbb{R}}\;\mathbb{C}\rightarrow {\mathbf T}_{\mathbb{C}} $ is an
isomorphism, so that there is a short exact sequence
$$
0\rightarrow \mathcal{K}_{k}\rightarrow \mathbf{T}_{k}
\stackrel{t}{\rightarrow }C_{k}(S^1)\rightarrow 0.
$$
where $\mathcal{K}_k$ is the $C^*$-algebra of compact operators on
a separable Hilbert space over $k$.

Let $\tau :C_k(S^1)\rightarrow k$ be the homomorphism given by $f\mapsto
f(1)$. Denote by $\mho_k$ the kernel of $\tau$. It is clear that
$\mho_{\mathbb{R}}$ is naturally isomorphic to the algebra
$C_0^{\mathbb{R}}(i\mathbb{R})$ defined in \cite{cunt}, and
$\mho_{\mathbb{C}}$ is isomorphic to $\Omega_{\mathbb{C}}$, where
$\Omega_k=\{f:I=[0,1]\rightarrow k \;|\;f(0)=f(1)=0\}$. The actions of $G$ on
the algebras considered above are trivial.

It follows from the Higson's theorem (see \cite{hig2}, and \cite{inka} for
the real case) that the functor $H_n$ is homotopy invariant for any $n\in
\mathbb{Z}$. Now, the proofs of the following proposition and theorem
coincide (up to trivial changes) with the proofs of suitable results in
\cite{cunt}.
\begin{proposition}
\label{ttob}
Let $\mathbb{H}$ be a stable cohomology theory and let
$g:\mathbf{T}_{k}\rightarrow k$ be the homomorphism defined by $v\mapsto 1$.
Then the homomorphism
$$
H_n(\id_A\otimes g):H_n(A)\xrightarrow{\simeq } H_n(A\otimes \mathbf{T}_{k})
$$
is an isomorphism, for any $A\in S_G$ and $n\in\mathbb Z$.
\end{proposition}

\begin{theorem}
\label{ztnuc} Let $\mathbb{H}$ be a stable cohomology theory,
$\mho A=A\otimes \mho _k$ and $\Omega A=A\otimes \Omega _k$. Then
there are natural isomorphisms
\begin{equation}
\label{isocunt} H_{n+1}(A)\simeq H_{n}(\mho A)
\;\;\;\text{and}\;\;\;H_{n-1}(A)\simeq H_n(\Omega A).
\end{equation}
for any $A\in S_G$ and $n\in \mathbb{Z}$.
\end{theorem}

As a consequence we have the following working principle:
\textit{two stable cohomology theories are isomorphic if and only
if they are isomorphic in some fixed dimension}. Thus in the
remaining sections of the paper we show that the families of
functors in Theorem \ref{corollary} are stable cohomology theories
and they are isomorphic when $n=0$.

In more detail, in section 2 we study an interpretation of algebraic and
topological $K$-theories of $C^*$-categories. Our definition is an adaptation
to our cases of some arguments from \cite{bass}, \cite{hig1}, \cite{quil}.
Let $A$ be a $C^*$-category and let $I$ be a $C^*$-ideal in $A$; let $a$ and
$a'$ be objects in $A$. We write $a\leq a'$ if there exists a morphism
$s:a\rightarrow a'$ such that $s^{*}s=1_a$ (such a morphism is said to be an
\textit{isometry}). Denote by $\mathcal{L}(a)$ (resp. by $\mathcal{I}(a)$)
the $C^{*}$-algebra $\hom_A(a,a)$ (resp $\hom_I(a,a)$). We have a
well-defined inductive system of abelian groups
$\{\mathbb{K}_{n}^a(\mathcal{L}(a)),\sigma _{aa'}\}_a$ and
$\{\mathbb{K}_{n}^t(\mathcal{I}(a)),\sigma _{aa'}\}_a$. We suppose that
\begin{equation}
\label{auxdef}
\mathbf{K}_{n}^a(A)=\underrightarrow{\lim}_a\mathbb{K}_{n}^a(
\mathcal{L}(a))\;\;\; \text{and}\;\;\;
\mathbf{K}_{n}^a(I)=\underrightarrow{\lim}_a\mathbb{K}_{n}^a(
\mathcal{I}(a))).
\end{equation}
Thanks to the results of A. Suslin and M. Wodzicki on the excision property
of algebraic $K$-groups on $C^*$-algebras \cite{suw}, the right hand side of
the second equation is well-defined. Algebraic $K$-groups obtained in this
way are naturally isomorphic to Quillen's $K$-groups $\mathbb{K}_{n}^Q(A)$
when $n\geq 0$; and are isomorphic to the Pedersen-Weibel $K$-groups in
negative dimensions. Note that a new interpretation of algebraic $K$-groups
implies existence of a simple flexible technical tool. Namely, any element of
an algebraic $K$-group of an additive category may be represented as an
element of an algebraic $K$-group of the endomorphism algebra of an object,
and such an interpretation is unique up to a manageable equivalence.
Throughout the paper, this principle will be used repeatedly.

In section 2, according to the excision property of algebraic $K$-groups on
the category of $C^*$-algebras \cite{suw}, we establish the excision property
for a short exact sequence associated to a $C^*$-ideal in an additive
$C^*$-category (see Proposition \ref{utul}). In section 3 this property is
used to prove Theorem \ref{exalg} about the excision property of functors
\begin{equation}
\label{katia} \{\mathbf{K}^a_n(\mathrm{Rep}(-;B))\}_{n\in Z}.
\end{equation}
In addition to the excision property, proof of Theorem \ref{exalg}
uses two nontrivial results. These are theorem \ref{topic} and
Theorem \ref{nichbisi}.

In section 4, the stability property of the functors \ref{katia}
will be shown.

Now, since the family of Kasparov's functors $KK_n(-;B),\;n\in\mathbb Z$ is a
stable cohomology theory \cite{cusk}, the proof of Theorem \ref{corollary}
boils down to showing the isomorphism
$$
\mathbb{K}_0(\mathrm{Rep}(A,B))\simeq KK_1(A,B),
$$
which is done in section 5.

\begin{remark}
\label{kark} Discussions for topological K-groups are omitted, because they
literally coincide with the considered case of algebraic $K$-groups.
\end{remark}

\begin{remark}
\label{kakan} The main result (with essential changes in definitions and
theorems) is also true for a locally compact group $G$. We hope to discuss
this case independently, not in this paper.
\end{remark}

\section{Some remarks on the algebraic $K$-theory of additive $C^*$-categories}

We will use the Pedersen-Weibel interpretation of algebraic $K$-groups
\cite{pwe}, denoted here by $\mathbb{K}_n^a$,  instead of Quillen's
definition of algebraic $K$-groups in \cite{quil} which was given through
homotopy groups of certain space. In this section we review some properties
of algebraic $K$-groups of idempotent-complete additive categories, based on
Pedersen-Weibel's nonconnective spectra (in this context there are defined
negative $K$-groups, too). Then, we reinterpret algebraic $K$-groups of
idempotent-complete additive $C^*$-categories and, with the aid of results
from \cite{suw}, generalize them to $C^*$-ideals in additive
$C^*$-categories. This material plays auxiliary role in this paper.

In the following lemma we list some simple properties of algebraic $K$-groups
which suffice for our purposes.

\begin{lemma}
\label{auxlem} Let $A$ be an idempotent-complete small additive category.
Then
\begin{enumerate}
\item
if $A=A_1\times A_2$, then
$\mathbb{K}^a_n(A)=\mathbb{K}^a_n(A_1)\times\mathbb{K}^a_n(A_2),\;\;
n\in\mathbb Z$;
\item
if $\{A_{\alpha }\}$ is a direct system of full additive
subcategories in $A$ such that $\underrightarrow{\lim}A_{\alpha}=A$.
Then
$\mathbb{K}^a_n(A)=\underrightarrow{\lim}\;\mathbb{K}^a_n(A_{\alpha}),\;\;n\in
\mathbb Z$.
\end{enumerate}
\end{lemma}

In the remaining part of this section we give an interpretation of the groups
$\mathbb{K}^a_n(A)$, $n\in \mathbb{Z}$ for additive $C^*$-categories. This
interpretation will convenient for our purposes.

Let $a$ be an object in an idempotent-complete $C^*$-category $A$. Consider
the full sub-$C^{*}$-category $A_a$ in $A$ consisting of all those objects
$a'$ of $A$ which admit an isometry $s:a'\rightarrow a^{\oplus_n}$. It is
clear that $A_a$ is an idempotent-complete $C^*$-category equivalent to the
category $\mathcal{P}(\mathcal{L}(a))$ of f.~g. projective modules over
$\mathcal{L}(a)$.

Now, consider a direct system of abelian groups for $A$ based on
subcategories $A_a$. It is evident that if there exists an isometry
$s:a'\rightarrow a$ then one has a natural additive inclusion $*$-functor
(not depending on $s$) $i_{a'a}:A_{a'}\rightarrow A_a$ and thus we have the
direct system $\{A_a,i_{a'a}\}_{(obA,\leq )}$ of idempotent-complete
$C^{*}$-categories. Because of the continuity property of algebraic
$K$-groups (property (2) in Lemma \ref{auxlem}) and the isomorphism of
categories $A=\underrightarrow{\lim}A_a$, one has an isomorphism
$$
\mathbb{K}^a_n(A)=\underrightarrow{\lim}\;\mathbb{K}^a_n(A_a)
$$
$n\in\mathbb Z$.

This suggests that $\mathbb{K}^a_n(A)$ can be interpreted in the
form
\begin{equation}
\label{kiso}
\mathbb{K}^a_n(A)=\underrightarrow{\lim}\;\mathbb{K}^a_n(\mathcal{L}(a)).
\end{equation}
Below it is done in detail.

\subsection{Algebraic $K$-functors of $C^{*}$-ideals \label{ctnuf}}
Let us make some comments on the results of A. Suslin and M. Wodzicki in
algebraic $K$-theory before we introduce our view on algebraic $K$-theory of
$C^{*}$-ideals. One of their main results is, by Proposition 10.2 in
\cite{suw}, that $C^*$-algebras have the factorization property
$(\mathbf{TF})_{\mathrm{right}}$. Thus any $C^{*}$-algebra possesses the
property $\mathbf{AH_{Z}}$. These results have many useful consequences in
algebraic $K$-theory of $C^*$-algebras, which are listed below. (Recall that
if $A$ is a $C^{*}$-algebra and $A^{+}$ is the $C^{*}$-algebra obtained by
adjoining a unit to $A$, then
\begin{equation}\label{Kcstar}
K_n^a(A)=\ker(K_n^a(A^+)\rightarrow K_n^a(k)),
\end{equation}
$n\in \mathbb{Z}$).

\begin{enumerate}
\item
$K^a_{i}$ is a covariant functor from the category of $C^{*}$-algebras and
$*$-homomorphisms to the category of abelian groups for any $i\in\mathbb{Z}$;
\item
For every unital $C^{*}$-algebra $R$ containing the $C^{*}$ algebra $A$ as a
two-sided ideal, the canonical map $K^a_n(A)\rightarrow K^a_{n}(R,A)$ is an
isomorphism;
\item
The natural embedding into the upper left corner $A\hookrightarrow M_{k}(A)$
induces, for every natural $n$, an isomorphism $K^a_{n}(A)\simeq
K^a_{n}(M_{k}(A))$;
\item
Any extension of $C^{*}$-algebras
$$
0\rightarrow I\rightarrow B\rightarrow A\rightarrow 0
$$
induces a functorial two-sided long exact sequence of algebraic $K$-groups
\begin{equation}\label{lsakt}
\cdots\rightarrow K^a_{i+1}(A)\rightarrow K^a_{i}(I)\rightarrow K^a_{i}(B)
\rightarrow K^a_{i}(A)\rightarrow \cdots  \;\;\;\;\;\; (i\in \mathbb{Z}).
\end{equation}
\item
Let $A$ be a $C^{*}$-algebra and let $u$ be a unitary element in a unital
$C^{*}$-algebra containing $A$ as a closed two-sided ideal. Then the inner
automorphism $ad(u):A\rightarrow A$ induces the identity map of algebraic
$K$-groups.
\end{enumerate}

Below we define algebraic $K$-groups for $C^*$-ideals. These groups possess
all properties similar to those represented above.

Let $A$ be an additive $C^{*}$-category and let $J$ be its closed
$C^{*}$-ideal. Let $\mathcal{L}_{A}(a)=\hom_A(a,a)$ and
$\mathcal{L}_{A}(a,J)=\mathcal{L}(a)_A\cap J$ for any object $a\in A$. The
latter is a closed ideal in the $C^{*}$-algebra $\mathcal{L}_{A}(a)$. Let us
write $a\leq a'$ if there is an isometry $v:a\rightarrow a'$, i.e
$v^{*}v=1_a$. The relation ``$a\leq a$'' makes the set of objects of $A$ into
a directed system. Any isometry $v:a\rightarrow a'$ defines a
$*$-homomorphism of $C^{*}$-algebras
$$
\mathrm{Ad}(v):\mathcal{L}_{A}(a)\rightarrow \mathcal{L}_{A}(a')
$$
by the rule $x\mapsto vxv^{*}$. It maps $\mathcal{L}_{A}(a,J)$ to
$\mathcal{L}_{A}(a',J)$.

Let $v_1:a\rightarrow a'$ and $v_2:a\rightarrow a'$ be two isometries. Then
$\mathrm{Ad}v_1$ and $\mathrm{Ad}v_2$ induce the same homomorphisms
$$
\mathrm{Ad}_{*}v_1=\mathrm{Ad}_{*}v_2:K_n^a(\mathcal{L}_{A}(a))\rightarrow
K_n^a(\mathcal{L}_{A}(a'))
$$
and
$$
\mathrm{Ad}_{*}v_1=\mathrm{Ad}_{*}v_2:K_n^a(\mathcal{L}_{A}(a,J))\rightarrow
K_n^a(\mathcal{L}_{A}(a',J)).
$$
The similar result for topological $K$-theory is in \cite{hig1}. This means
that the homomorphism $\nu_{*}^{aa'}=K_n^a(\nu^{aa'})$ is independent of
choosing an isometry $\nu^{aa'}:a\rightarrow a'$. Therefore one has a
directed system $\{K_n^a(\mathcal{L}_{A}(a,J)),\nu_{*}^{aa'})\}_{a,a'\in
obA}$ of abelian groups, for all $n\in\mathbb Z$.

\begin{definition}
\label{wb} Let $A$ be an additive $C^*$-category and let $J$ be its closed
$C^{*}$-ideal. Define
\begin{equation}
\label{ktg}
\begin{array}{c}
\mathbf{K}_n^a(A,J)=\underrightarrow{\lim
}K_n^a(\mathcal{L}_{A}(a,J))\;\;\;\text{and}\;\;\;\mathbf{K}_n^a(J)=\mathbf{K}_n^a(M(J),J)
\end{array}
\end{equation}
where $M(J)$ is  the so called \textit{multiplier} $C^*$-category of
$J$ \cite{kan3}.
\end{definition}

\begin{lemma}
\label{wi} Let $J$ be a $C^{*}$-ideal in an additive
$C^{*}$-category $A$. Then
\begin{enumerate}
\item
$\mathbf{K}_n^a(J)=\mathbf{K}_n^a(A,J)$;
\item
if $A'$ is a cofinal subcategory in $A$ then $\mathbf{K}_n^a(A')\simeq
\mathbf{K}_n^a(A)$;
\item
if $A$ is idempotent-complete, then
$\mathbf{K}_n^a(A)\simeq\mathbb{K}_n^a(A)$.
\end{enumerate}
\end{lemma}

\begin{proof}
1. The natural $*$-functor $\rho:A\rightarrow M(J)$ induced from the identity
on $J$ by the universal property of $M(J)$ obviously preserves the relation
``$\leq$''. This implies that there is a natural morphism of directed systems
$$
\{K_n^a(\mathcal{L}_{A}(a,J)),\nu_{*}^{aa'})\}_{a,a'\in\mathrm{ob}
A}\xrightarrow{\{\rho^a_n\}}
\{K_n^a(\mathcal{L}_{M(J)}(a,J)),\nu_{*}^{aa'})\}_{a,a'\in\mathrm{ob} M(J)},
$$
where the homomorphism $\rho^a_n:K_n^a(\mathcal{L}_{A}(a,J))\rightarrow
K_n^a(\mathcal{L}_{M(J)}(a,J))$ is induced by the above $*$-functor
$\rho:A\rightarrow M(J)$. In view of the isomorphism $K_n(A)\rightarrow
K_{n}(R,A)$ \cite{suw} one concludes that the homomorphism $\rho^a_n$ is an
isomorphism for all $n\in\mathbb{Z}\;\;a\in A$. This morphism is cofinal,
since if ``$a\leq a'$'' in $M(J)$ then ``$a\leq a\oplus a'$'' in $A$ and
``$a'\leq a\oplus a'$'' in $M(J)$. Therefore, $\{\rho^a_n\}$ is an
isomorphism of direct systems.

2. This is an easy consequence of Definition \ref{wb}.

3. This results from comparison of Definition \ref{wb} and
isomorphism \ref{kiso}.
\end{proof}

Now, we prove the excision property of algebraic $K$-theory which will be
used in the next section.
\begin{proposition}
\label{utul} Let $A$ be an additive $C^{*}$-category and let $J$ be a
$C^*$-ideal in $A$. Then the two-sided sequence of algebraic $K$-groups
\begin{equation}
\label{exact} ...\rightarrow \mathbf{K}_{n+1}^a(A/J)\rightarrow
\mathbf{K}_{n}^a(J)\rightarrow \mathbf{K}_{n}^a(A)\rightarrow
\mathbf{K}_n^a(A/J)\rightarrow ...,\;\;\;
\end{equation}
$n\in \mathbb Z$, is exact.
\end{proposition}

\begin{proof}
Consider the exact sequence of $C^{*}$-algebras
$$
0\rightarrow \mathcal{L}(a,J)\rightarrow
\mathcal{L}(a,A)\rightarrow
\mathcal{L}(a,A)/\mathcal{L}(a,J)\rightarrow 0.
$$
By the excision property of algebraic $K$-theory on $C^{*}$-algebras
\cite{suw}, one has a two-sided long exact sequence of algebraic $K$-groups
\begin{multline}\label{exb} \cdots \rightarrow
K_n^a(\mathcal{L}(a,A)/\mathcal{L}(a,J))\rightarrow \\
K_{n-1}^a(\mathcal{L}(a,J))\rightarrow
K_{n-1}^a(\mathcal{L}(a,A))\rightarrow K_{n-1}^a(
\mathcal{L}(a,A)/\mathcal{L}(a,J))\rightarrow \cdots
\end{multline}
According to the exact sequence \ref{exb} and the fact that directed colimits
preserve exactness, one obtains the following long exact sequence of abelian
groups:
\begin{multline} \cdots \rightarrow
\underrightarrow{\lim}K_{n+1}^a(\mathcal{L}(a,B)/\mathcal{L}(a,J))\rightarrow \\
\mathbf{K}_{n}^a(J)\rightarrow \mathbf{K}_{n}^a(B)\rightarrow
\underrightarrow{\lim}K_{n}^a(
\mathcal{L}(a,B)/\mathcal{L}(a,J))\rightarrow \cdots
\end{multline}
There is a natural morphism of directed systems
$$
\{\omega _a\}:\{K_{n}^a(
\mathcal{L}(a,A)/\mathcal{L}(a,J)\}\rightarrow \{K_{n}^a(
\mathcal{L}(a,A/J)\}
$$
so that $\omega_a$ is the identity map for any object $a$ in $A$.
It is clear that this morphism is cofinal. Thus the induced
homomorphism
$$
\omega:\underrightarrow{\lim}K_{n}^a(
\mathcal{L}(a,A)/\mathcal{L}(a,J))\rightarrow
\underrightarrow{\lim}K_{n}^a(
\mathcal{L}(a,A/J))=\mathbf{K}_{n}^a(A/J)
$$
is an isomorphism.
\end{proof}

\section{Excision Property of $\mathbb{K}^a_n((Rep(-;B))$ \label{auxiso}}

In this section we will prove that the contravariant functors
$\mathbb{K}^a_n((\mathrm{Rep}(-;B))$, $n\in \mathbb{Z}$, have the weak
excision property. The similar result for topological $K$-theory, in a
particular case, has been proved in \cite{hig1}.

\begin{remark}
\label{wepe} From now on for convenience of calculations the
functors $\mathbf{K}^a_n((Rep(-;B))$ are considered instead of
$\mathbf{K}^a_n((\mathrm{Rep}(-;B))$ and
$\mathbb{K}^a_n((\mathrm{Rep}(-;B))$. Since $Rep(A;B)$ is a
cofinal full subcategory in $\mathrm{Rep}(A;B)$, by Lemma \ref{wi}
(2),(3) these functors are isomorphic.
\end{remark}

For any closed invariant ideal $J$ in a separable $C^*$-algebra $A$ from
$S_G$ there is a $C^*$-ideal $D(A,J;B)$ in $Rep(A;B)$ which is defined in the
following manner. Let $(E,\phi)$ and $(E',\phi')$ be objects in $Rep(A,B)$. A
morphism $\alpha:(E,\phi)\rightarrow (E',\phi ')$ in $Rep(A,B)$ is in
$D(A,J;B)$ if
$$
\alpha \phi (x)\in \mathcal{K}((E,\phi ),(E'\phi
'))\;\;\text{and}\;\;\phi '(x)\alpha \in \mathcal{K}((E,\phi
)),\;\text{for all}\;x\in A.
$$
The space of all morphisms from $(E,\phi)$ to $(E',\phi')$ in the $C^*$-ideal
$D(A,J;B)$ will be denoted by $D_{\phi,\phi'}(A,J;E,E';B)$ (if
$(E',\phi')=(E,\phi)$ then it is also denoted by $D_\phi(A,J;E;B)$) (cf.
\cite{hig1}) .

\begin{theorem}
\label{excision}Let $B$ be a $\sigma$-unital $C^{*}$-algebra and let
$0\rightarrow I\rightarrow A\stackrel{p}{\rightarrow }A/I\rightarrow 0$ be a
proper sequence of separable $C^{*}$-algebras in $S_G$. Then the sequence of
groups
\begin{multline}
\label{exalg} _...\rightarrow
\mathbf{K}_{n}^a(Rep(A/J,B))\rightarrow
\mathbf{K}_n^a(Rep(A,B))\rightarrow
 \mathbf{K}_n^a(Rep(J,B))\xrightarrow{\partial} \\
 \xrightarrow{\partial}\mathbf{K} _{n-1}^a(Rep(A,B))\rightarrow ...
\end{multline}
is exact, for all $n\in\mathbb Z$.
\end{theorem}

\begin{proof}
Consider the short exact sequence of $C^{*}$-categories and of a $C^*$-ideal
$$
0\rightarrow D(A,J;B)\rightarrow Rep(A,B)\rightarrow
Rep(A,B)/D(A,J;B)\rightarrow 0.
$$
According to Proposition \ref{utul}, one has a two-sided long
exact sequence
\begin{multline}\label{esk}
...\rightarrow \mathbf{K} _{n}^a(D(A,J;B))\rightarrow
\mathbf{K}_n^a(Rep(A,B))\rightarrow \\\rightarrow \mathbf{K}
_n^a(Rep(A,B)/D(A,J;B))\stackrel{\partial }{\rightarrow }
\mathbf{K} _{n-1}^a(D(A,J;B))\rightarrow...
\end{multline}
of abelian groups. According to Theorems \ref{topic} and
\ref{nichbisi}, replacements
\begin{equation*}
\mathbf{K}_n^a(Rep(A;B)/D(A,J;B))\;\; \text{by}\;\;\;
\mathbf{K}_n^a(Rep(J;B))
\end{equation*}
and
\begin{equation*}
\mathbf{K}_n^a(D(A,J;B)\;\;\text{by}\;\;\;\mathbf{K}_n^a(Rep(A/J;B)),
\end{equation*}
ensure exactness of the two-sided long exact sequence (\ref {exalg}). Thus it
suffices to prove Theorems \ref{topic} and \ref{nichbisi}, which will be done
in the next part of this section.
\end{proof}

\subsection{On the Isomorphism
$\mathbf{K}_n^a(Rep(A;B)/D(A,J;B))\approx
\mathbf{K}_n^a(Rep(J;B))$}

Let $(E,\phi)$ be an object in $Rep(A,B)$ and let $j:J\rightarrow A$ be the
natural equivariant inclusion. There is a $*$-functor induced by the natural
inclusion $j$
\begin{equation}\label{xxx}
{\mathrm j}:Rep(A;B)\rightarrow Rep(J;B)
\end{equation}
defined by the assignments $(E,\phi)\mapsto (E,\phi j)$ and $x\mapsto x$.

The following trivial lemma is used in the next proposition.
\begin{lemma}
\label{print} Let $A$ and $B$ be additive $C^*$-categories and let
$F:A\rightarrow B$ be an additive $*$-functor. Then $F$ is a
$*$-isomorphism if and only if $F$ is bijective on objects and the
induced $*$-homomorphisms of $C^*$-algebras
$F_a:\mathcal{L}(a)\rightarrow \mathcal{L}(f(a))$ are
$*$-isomorphisms for all objects $a$ in $A$.
\end{lemma}
The following proposition is a slight generalization of the similar result in
\cite{hig1}).
\begin{proposition}
\label{tipo} The canonical $*$-functor \ref{xxx} maps $D(A,J;B)$ to
$D(J,J;B)$ and the induced $*$-functor
\begin{equation}\label{yxy}
\xi:Rep(A;B)/D(A,J;B)\rightarrow Rep(J;B)/D(J,J;B)
\end{equation}
is an isomorphism of $C^{*}$-categories.
\end{proposition}

\begin{proof} (cf. \cite{hig1})
By lemma \ref{print} it suffices to show that for any object $(E,\phi)$ the
$*$-homomorphism of $C^*$-algebras
$$
\xi_{J,\phi}:D_\phi(A;E;B)/D_\phi (A,J,E;B)\rightarrow D_{\phi\cdot
j}(J,E;B)/D_{\phi\cdot j}(J,J,E;B)
$$
is a $*$-isomorphism. It is easy to show that $\xi_{J,\phi}$ is a
monomorphism. To show that $\xi_{J,\phi}$ is an epimorphism, take $x\in
D_{\phi\cdot j}(J,E;B)$ and let $E_1$ be the $G$-$C^{*}$-algebra in
$\mathcal{L}(E)$ generated by $\phi(J)\cup\mathcal{K}(E)$; let $E_2$ be the
separable $G$-$C^{*}$-algebra generated by all elements of the form $[x,\phi
(y)]$, $y\in J$; and let $\mathcal{F}$ be the $G$-invariant separable linear
space generated by $x$ and $\phi(A)$. One has
\begin{itemize}
\item  {$E_1\cdot E_2\subset \mathcal{K}(E)$,\ \ \ \ because $\phi
(b)[\phi (a),x]\sim [\phi (ba),x]\in \mathcal{K}(E),\;\;$ $a\in
A,\;b\in J$, }

\item  {$[\mathcal{F},E_1]\subset E_1$, \ \ \ \ \ \ because
$[x,\phi (J)]\subset \mathcal{K}(E)$ and $[\phi (A),\phi
(J)]\subset \phi (J)$.}
\end{itemize}

From a technical theorem by Kasparov it follows that there exists a positive
$G$-invariant operator $X$ such that
\begin{enumerate}
\item
$X\cdot\phi(J)\subset\mathcal{K}(E)$;
\item
$(1-X)\cdot[\phi(A),x]\subset\mathcal{K}(E)$;
\item
$[x,X]\in\mathcal{K}(E)$.
\end{enumerate}

Since $[(1-X)x,\phi(a)]=(1-X)[x,\phi (a)]-[X,\phi(a)]x$, it follows from (2)
and (3) that $(1-X)x\in D_\phi(A,E;B)$. In addition, it follows from (2) that
$Xx\in D_{\phi\cdot j}(J,J,E;B)$, and so that the image of $(1-X)x$ in
$D_{\phi\cdot j}(J,E;B)/D_{\phi\cdot j}(J,J,E;B)$ coincides with the image of
$x$.
\end{proof}

Now, we prove the following

\begin{theorem}
\label{topic} Let $A$ be a separable $C^{*}$-algebra and let $B$ be a $\sigma
$-unital $C^{*}$-algebra. Let $J$ be a closed ideal in $A$. There exists an
essential isomorphism
\begin{equation}
\label{qjp} {\bf K}_n^a(Rep(A,B)/D(A,J;B))\approx{\bf K}_n^a(Rep(J,B))
\end{equation}
\end{theorem}

\begin{proof}
According to Proposition \ref{tipo}, it suffices to show that the
homomorphism
$$
\mathbf{K}_{*}^a(Rep(J;B))\rightarrow
\mathbf{K}_{*}^a(Rep(J;B)/D(J,J;B))
$$
is an isomorphism. The exact sequence
\begin{multline}\label{eska}
...\rightarrow \mathbf{K} _{n}^a(D(J,J;B))\rightarrow
\mathbf{K}_n^a(Rep(J,B))\rightarrow \\\rightarrow \mathbf{K}
_n^a(Rep(J,B)/D(J,J;B))\stackrel{\partial }{\rightarrow }
\mathbf{K} _{n-1}^a(D(J,J;B))\rightarrow...
\end{multline}
shows that it suffices to show $\mathbf{K}^a_{*}(D(J,J;B))=0$. According to
Kasparov's stabilization theorem, one concludes that the full subcategory
$Rep_{H_B}(J;B)$ on all objects of the form $(H_B,\varphi)$ is a cofinal
subcategory in $Rep(J;B)$, where $H_B$ is Kasparov's universal Hilbert
$B$-module \cite{kas1}. Note that the canonical isometry
$$
i_1:H_B\rightarrow H_B\oplus H_B
$$
in the first summand is in $D_{\phi ,\phi \oplus
0}(J;H_B,H_B\oplus H_B;B)$. So it induces inner homomorphism
$$
\mathrm{ad}(i_1):D_\phi (J,J;H_B;B)\rightarrow D_{\phi \oplus
0}(J,J;H_B\oplus H_B;B).
$$
Consider the sequence of $*$-homomorphisms
\begin{equation}
D_\phi(J,J;H_B;B)\rightarrow D_{\phi\oplus\phi}(J,J;H_B\oplus H_B;B)\subset
D_{\phi\oplus0}(J,J;H_B\oplus H_B;B),
\end{equation}
where the inclusion is given by $x\mapsto x$. If the first arrow is induced
by the inclusion $\iota_1:H_B\rightarrow H_B\oplus H_B$ into the first
summand, then the composite is $\mathrm{ad}(i_1)$. If the first arrow is
induced by the inclusion $\iota_2:H_B\rightarrow H_B\oplus H_B$ into the
second summand, one obtains a homomorphism $\lambda$. Since
$K_n^a(\mathrm{ad}(\iota_1))=K_n^a(\mathrm{ad}(\iota_2))$, one has
$K_n^a(\mathrm{ad}(i_1))=K_n^a(\lambda)$. On the other hand, the homomorphism
$\lambda$ is the composite of $*$-homomorphisms of $C^*$-algebras
$$
D_\phi(J,J;H_B;B)\rightarrow D_0(J,J;H_B;B)\rightarrow
D_{\phi\oplus0}(J,J;H_B\oplus H_B;B),
$$
defined by the assignments
$$
x\mapsto x \textrm{ and } x\mapsto\left(
\begin{array}{cc}
0 & 0 \\
0 & x
\end{array}
\right).
$$
Note that one has $D_0(J,J;H_B;B)\simeq M(J\otimes \mathcal{K}_G)$. It is a
well-known fact that the latter algebra has trivial algebraic $K$-groups. If
we apply $K$-functors then the homomorphism corresponding to $\lambda$ will
be zero. Now, as a consequence we have the following. Let $\alpha\in
K_{n}^a(D_\phi(J,J;H_B;B))$ represent an element in
$\mathbf{K}_n^a(D(J,J;B))$. Since
$$
\alpha=K_{n}^a(\mathrm{ad}(i_1)(\alpha ))=K_{n}^a (\lambda (\alpha))=0,
$$
one concludes that the class of $\alpha$ in $\mathbf{K}_{*}^a(D(J,J;B))$ is
zero. Therefore $\mathbf{K}_{*}^a(D(J,J;B))=0$.
\end{proof}

\subsection{On the Isomorphism $\mathbf{K}_n^a(Rep(A/J;B))\simeq
\mathbf{K}_n^a(D(A,J;B))$}

Let
$$
0\rightarrow J\xrightarrow{j} A\xrightarrow{q} A/J\rightarrow 0
$$
be a proper exact sequence and let $\sigma:A/J\rightarrow A$ be a completely
positive and contractive (equivariant) section.

Let $(E,\phi)$ be an object in $Rep(A;B)$. A $*$-homomorphism
$$
\psi =\left(
\begin{array}{cc}
\psi _{11} & \psi _{12} \\
\psi _{21} & \psi _{22}
\end{array}
\right) :A/J\rightarrow \mathcal{L}(E\oplus E'))
$$
will be called a \emph{$\sigma$-dilation} for $\phi$ if
$\psi_{11}(a)=\phi(\sigma(a))$, where $E'$ is a right Hilbert $B$-module. By
generalized Stinespring's theorem there exists a $\sigma$-dilation for
$\phi$, where $\sigma:A/J\rightarrow A$ is a completely positive and
contractive section \cite{kas1}.

\begin{lemma}
\label{zb}Let $\psi$ be a $\sigma$-dilation for $\phi$. Then
\begin{enumerate}
\item
$\psi_{12}(a^{*})=\psi_{21}(a)^{*}$;
\item
for any $a,b\in A/J$ there exists a $j\in J$ such that
$\psi_{12}(a)\psi_{21}(b)=\phi(j)$;
\item
$\psi_{12}(a)x$ and $x\psi_{21}(a)$ are compact morphisms for any $a\in A/J$
and $x\in D_\phi(A,J;B)$;
\item
Let $\phi:A/J\rightarrow\mathcal{L}(E)$ be a $*$-homomorphism and let
$\psi:A/J\rightarrow \mathcal{L}(E\oplus E')$ be a $\sigma$-dilation for
$\phi q$. There exists a $*$-homomorphism $\varphi:A/J\rightarrow
\mathcal{L}(E')$ such that
\end{enumerate}
$$
\psi=
\begin{pmatrix}
\phi & 0\\
0 & \varphi
\end{pmatrix}.
$$
\end{lemma}

\begin{proof}
\textit{The case (1)} is trivial, because $\psi $ is a $*$-homomorphism.

\textit{The case (2)}. Since $\psi$ is $\sigma$-dilation for $\phi$, one has
$\phi\cdot(\sigma(ab)-\sigma(a)\cdot\sigma(b))=\psi_{12}(a)\cdot\psi_{21}(b)$.
But $j=\sigma(ab)-\sigma(a)\cdot s(b)\in J$. Therefore $\psi_{12}(a)\cdot
\psi_{21}(b)=\phi(j)$.

\textit{The case (3)}. If $x\in D_\phi(A,J;B)$ then, by definition, $x\phi
(j)$ and $\phi(j)x$ are compact morphisms for any $j\in J$. According to
cases (1) and (2), one has
$x\psi_{12}(a)\cdot\psi_{12}^{*}(a)x^{*}=x\phi(j')x^{*}$\ for some $j'\in J$.
Thus $x\psi_{12}(a)\cdot\psi_{12}(a^*)x^*$ is a compact morphism. Therefore
$x\psi_{12}(a)$ and $\psi_{21}(a)x$ ($=(x^*\psi_{12}(a^*))^*$) are compact
morphisms, too.

\textit{The case (4)}. Since $\psi$ is a $\sigma$-dilation for $\phi q$,
$$
\psi=
\begin{pmatrix}
\phi &\psi_{12}\\
\psi_{21} & \varphi
\end{pmatrix}.
$$
According to case (2), for any $a,b\in A/J$ there exists $j\in J$ such that
$\psi_{12}(a)\psi_{21}(b)=\phi q(j)=0$. Applying the case (1), one has
$\psi_{12}(a)\psi_{12}^*(a)=0$. Therefore $\psi_{12}(a)=0$ and
$\psi_{21}(a)=0$ for all $a\in A/J$ and $\varphi$ is a $*$-homomorphism.
\end{proof}

The following lemma is used in Theorem \ref{nichbisi}.
\begin{lemma}
\label{zs} Let $(E,\phi)$ be an object in $Rep(A,B)$ and let
$\psi:A/J\rightarrow\mathcal{L}(E\oplus E')$ be a $\sigma$-dilation for
$\phi$. Then a map $x=\left(
\begin{array}{cc}
x_{11} & x_{12} \\
x_{21} & x_{22}
\end{array}
\right) \mapsto x'=\left(
\begin{array}{ccc}
x_{11}&0&x_{12}\\
0&0&0\\
x_{21}&0&x_{22}
\end{array}
\right)$ defines a $*$-monomorphism
\begin{equation}
\label{ksi}\rho:M_2(D_\phi(A,J,E\oplus E;B))\rightarrow D_{\psi\cdot q\oplus
\phi}(A,J,E\oplus E'\oplus E;B).
\end{equation}
\end{lemma}

\begin{proof}
It suffices to show that $x'\in D_{\psi\cdot q\oplus\phi}(A,J,E\oplus
E'\oplus E;B).$ By assumption one has
$$
(\phi(a)\oplus\phi(a))x-x(\phi(a)\oplus\phi(a))\in\mathcal{K}(E\oplus E),
$$
for any $a\in A$, and
$$
(\phi(b)\oplus\phi(b))x\in\mathcal{K}(E\oplus E),
$$
$x(\phi(b)\oplus\phi(b))\in\mathcal{K}(E\oplus E)$ for any $b\in J$. It
implies that
$$
\begin{array}{c}
\phi(a)x_{mn}-x_{mn}\phi(a)\in\mathcal{K}(E),\;\;\; \phi(b)x_{mn}\in
\mathcal{K}(E),\;\;\;x_{mn}\phi(b)\in\mathcal{K}(E),\;\;
\end{array}
$$
$a\in A$ and $b\in J$. Then $(\psi(q(a))\oplus \phi(a))\cdot
x'-x'\cdot(\psi(q(a))\oplus\phi(a))=$
$$
=\left(
\begin{array}{ccc}
x_{11}\psi_{11}(q(a))-\psi_{11}(q(a))x_{11}&x_{11}\psi_{12}(p(a))&
x_{12}\phi(a)-\psi_{11}(q(a))x_{12}\\
\psi_{21}(q(a))x_{11}&0&\psi_{21}(q(a))x_{12}\\
x_{21}\psi_{11}(q(a))-\phi(a)x_{21}&x_{21}\psi_{12}(q(a))&x_{22}\phi(a)-\phi(a)x_{22}
\end{array}
\right).
$$
By Lemma \ref{zb} (3), the morphisms $\psi_{21}(q(a))x_{11}$, $x_{11}\psi
_{12}(q(a))$, $x_{21}\psi _{12}(q(a))$ and $\psi _{21}(q(a))x_{12}$ are
compact. Using the fact that $\phi(a)-\psi_{11}(q(a))\in\phi(J)$, one has
$$
(\psi(q(a))\oplus\phi(a))\cdot x'-x'\cdot(\psi(p(a))\oplus\phi(a))\in
\mathcal{K}(E\oplus E'\oplus E),\;\;a\in A.
$$
To show that $(\psi(q(b))\oplus \phi(b))\cdot x'$ and
$x'\cdot(\psi(q(b))\oplus\phi(b))$ are in $\mathcal{K}(E\oplus E'\oplus E)$
when $b\in J$, note that $(\psi(q(b))\oplus\phi(b))\cdot x'$ and
$x'\cdot(\psi(q(b))\oplus\phi(b))$ are equal to
$$
\left(
\begin{array}{ccc}
0&0&0\\
0&0&0\\
\phi(b)x_{21}&0&\phi(b)x_{22}
\end{array}
\right) \;\;{\rm and}\;\;\left(
\begin{array}{ccc}
0&0&x_{12}\phi(b)\\
0&0&0\\
0&0&x_{22}\phi(b)
\end{array}
\right)
$$
respectively. They are compact morphisms because each entry of matrices is a
compact morphism.
\end{proof}

The category $Rep(A,B)$ may be identified with the full $C^{*}$-subcategory
$D^{(q)}(A,J;B)$ in $D(A,J;B)$, on the objects all pairs of the form
$(E,\phi\cdot q)$, considering $(E,\phi)$ as an object in $Rep(A/J;B)$. Let
$$
\varepsilon:D^{(q)}(A,J;B)\hookrightarrow D(A,J;B)
$$
be the natural inclusion.

\begin{theorem}
\label{nichbisi} Let $0\rightarrow J\xrightarrow{j}
A\xrightarrow{q} A/J\rightarrow 0$ be a proper exact sequence of
separable $C^{*}$-algebras. Then the induced homomorphism
$$
\Gamma_n:{\bf K}_n^a(Rep(A/J;B))\rightarrow {\bf K}_n^a(D(A,J;B))
$$
is an isomorphism for all $n\in\mathbb Z$.
\end{theorem}

\begin{proof}
According to the discussion above, it suffices to show that $\varepsilon$
induces an isomorphism
$$
\varepsilon_n:{\bf K}_n^a(D^{(q)}(A,J;B))\rightarrow
\mathbf{K}_n^a(D(A,J;B)),
$$
for all $n\in\mathbb Z$.

(1). $ \varepsilon_n $ \textit{is a monomorphism}.

Let $(E,\phi\cdot q)$ be an object in $D^{(q)}(A,J;B)$ and suppose that the
class of an element
$$
y\in K_n^a(D_{\phi \cdot q}(A,J;E,B))
$$
in $\mathbf{K}_n^a(D(A,J;B))$ is zero. This means that there exists
an isometry
$$
s:(E,\phi\cdot q)\rightarrow(E',\psi)
$$
in $Rep(A;B)$ such that Ad$(s)_n(y)=0$ in $K_n^a(D_{\psi}(A,J;E',B))$.

According to Lemma \ref{zb}, it is easy to show that
$$
s'=\left(%
\begin{array}{cc}
s&0\\
0&0
\end{array}%
\right):(E\oplus H_B,\eta\cdot q)\rightarrow(E'\oplus H_B,\eta'\cdot q)
$$
is an isometry in $Rep(A/J,B)$, where $\eta$ and $\eta'$ are
$\sigma$-dilations of $\phi\cdot q$ and $\psi$ respectively. By Lemma
\ref{zb} $\eta$ has form
$$
\eta=\left(
\begin{array}{cc}
\phi&0\\
0&\chi
\end{array}
\right),
$$
where $\chi$ is a $*$-homomorphism from $A$ to $\mathcal{L}(H_B)$. There is a
homomorphism
$$
\nu:D_{\psi}(A,J;E',B)\rightarrow D_{\eta'q}(A,J;E'\oplus\mathcal{H}_B,B),
$$
defined by $x\mapsto\left(
\begin{array}{cc}
x&0\\
0&0
\end{array}
\right)$, and an isometry
$$
i_1:(E,\phi)\rightarrow(E\oplus H_B,\eta=\left(
\begin{array}{cc}
\phi&0\\
0&\chi
\end{array}
\right))
$$
(inclusion into the first summand). It is clear that
Ad$(s'i_1)=\nu\mathrm{Ad}(s)$. Therefore Ad$(s'i_1)_n(y)=0$ in $D_{\eta
'q}(A,J;E'\oplus \mathcal{H}_B,B)$ and in $\mathbf{K}_n^a(D^{(q)}(A,J;B))$
too. This means that the class of $y$ in $\mathbf{K}_n^a(D^{(q)}(A,J;B))$ is
zero.

(2). $ \varepsilon_n$ \textit{is an epimorphism}.

Let an element in $\mathbf{K}_n^a(D(A,J;B))$ be represented by an element
$x\in K_n^a(D_{\phi }(A,J;E,B))$. Consider the $*$-homomorphism
$$
\theta:D_\phi(A,J;E;B)\rightarrow D_{\psi q}(A,J;E\oplus H_B,B)
$$
given by
$$
x\mapsto\left(
\begin{array}{cc}
x&0\\
0&0
\end{array}
\right),
$$
where $\psi$ is a $\sigma$-dilation of $\phi$.

Let us show that classes of elements $\theta_n(x)$ and $x$ in
$\mathbf{K}_n^a(D(A,J;B))$ coincide. Let $\vartheta$ be the composite
$$
D_\phi(A,J;E;B)\xrightarrow{\theta}D_{\psi q}(A,J;E\oplus
H_B;B)\xrightarrow{\textrm{Ad}(i_{E\oplus H_B})}D_{\psi
q\oplus\phi}(A,J;E\oplus H_B\oplus E;B),
$$
where the second arrow is induced by the isometry into the first two
summands. It is clear that the $*$-homomorphism $\vartheta$
coincides with the $*$-homomorphism defined by
$$
x\mapsto\left(
\begin{array}{ccc}
x&0&0\\
0&0&0\\
0&0&0
\end{array}
\right).
$$
On the other hand this homomorphism may be interpreted as the composite
$$
D_\phi(A,J;E,B)\xrightarrow{i_1}M_2(D_\phi(A,J;E,B))
\xrightarrow{\rho}D_{\psi\cdot q\oplus\phi}(A,J;E\oplus H_B\oplus E,B)
$$
where $i_1$ is given by $x\mapsto\left(
\begin{array}{cc}
x&0\\
0&0
\end{array}
\right)$ and $\rho$ is defined by
$$
x=\left(
\begin{array}{cc}
x_{11}&x_{12}\\
x_{21}&x_{22}
\end{array}
\right)\mapsto x'=\left(
\begin{array}{ccc}
x_{11}&0&x_{12}\\
0&0&0\\
x_{21}&0&x_{22}
\end{array}
\right),
$$
as in Lemma \ref{zs}. Consider another homomorphism $\eta$---the composite of
the sequence of $*$-homomorphisms
$$
D_\phi(A,J;E,B)\xrightarrow{i_2}M_2(D_\phi(A,J;E,B))
\xrightarrow{\rho}D_{\psi\cdot q\oplus\phi}(A,J;E\oplus H_B\oplus E,B)
$$
where $i_2$ is given by $x\mapsto\left(
\begin{array}{cc}
0&0\\
0&x
\end{array}
\right)$. Since $(i_1)_n=(i_2)_n$, $n\in\mathbb Z$, one has
$(\vartheta)_n=(\eta)_n$. But $(\eta)_n=(\mathrm{ad}(i_3))_n$, $n\in\mathbb
Z$, where $i_3:E\rightarrow E\oplus H_B\oplus E$ is an isometry in the third
summand. Therefore, classes of the elements $\vartheta_n(x)$,
$(\mathrm{ad}(i_3))_n(x)$, $x$ and $\theta_n(x)$ in
$\mathbf{K}_n^a(D(A,J;B))$ coincide. Therefore the class of the element
$\theta_n(x)$ in $\mathbf{K}_n^a(D^{(q)}(A,J;B))$ is the desired element.
\end{proof}

\section{Stability Property of $\mathbb{K}^a_n(\mathrm{Rep}(-;B))$ \label{stab}}

Everywhere below in this section, $\mathcal{K}$ is the $C^*$-algebra of
compact operators on a countable generated Hilbert space $\mathcal{H}$
considered as an object of $S_G$ via trivial action of the compact group $G$.

Let $p\in\mathcal{K}$ be a rank one projection and let $A$ be a $C^*$-algebra
in $S_G$; let $e_{A}:A\rightarrow A\otimes\mathcal{K}$ be the
$*$-homomorphism defined by $a\mapsto a\otimes p$, $a\in A$. Then one has the
induced functor
\begin{equation}
\label{stabKR}e^*_{A}:Rep(A\otimes\mathcal{K};B)\rightarrow Rep(A;B),
\end{equation}
defined by assignments $(E,\varphi)\mapsto(E,\varphi e_A)$ (on objects) and
$x\mapsto x$ (on morphisms).

There is a $*$-functor
$$
\varepsilon^A:Rep(A;B)\rightarrow Rep(A\otimes\mathcal{K};B)
$$
defined by assignments $(E,\phi)\mapsto(E\otimes _k\mathcal{H},\phi\otimes
\id_{\mathcal{K}})$ (on objects) and $f\mapsto f\otimes\id_{\mathcal{H}}$ (on
morphisms). Indeed, let $(E,\phi)$ be an object in $Rep(A;B)$. One has the
induced $*$-homomorphism
$$
\phi\otimes\id_{\mathcal{K}}:A\otimes\mathcal{K}\rightarrow
\mathcal{L}(E\otimes_k\mathcal{H}).
$$
Let $f:(E,\phi)\rightarrow(E',\phi')$ be a morphism in $Rep(A;B) $, i.~e.
$f\phi(a)-\phi'(a)f\in\mathcal{K}(E,E')$, $a\in A$. Then
\begin{multline}
(f\otimes\id_{\mathcal{H}})((\phi\otimes
\id_{\mathcal{K}})(a\otimes\kappa))-((\phi\otimes
\id_{\mathcal{K}})(a\otimes \kappa))(f\otimes\id_{\mathcal{H}})=\\
=(f\phi(a)-\phi (a)f)\otimes \kappa \in \mathcal{K}(E\otimes
_k\mathcal{H},E'\otimes _k\mathcal{H})
\end{multline}
for all $a\in A$, $\kappa\in\mathcal{K}$.

Now, in view of Remark \ref{wepe} stability property of
$\mathbb{K}^a_n(\mathrm{Rep}(-;B))$ may be formulated as follows.

\begin{theorem}
\label{sdfm} For any rank one projection $p\in\mathcal{K}$ and any
$C^*$-algebra $A$ in $S_G$ the homomorphism
$$
e_n^{A}=\mathbf{K}^a_n(e^*_{A}):\mathbf{K}^a_n(Rep(A\otimes
\mathcal{K};B))\rightarrow \mathbf{K}^a_n(Rep(A;B)),
$$
induced by the functor \ref{stabKR} is an isomorphism.
\end{theorem}

\begin{proof}
Let $\varepsilon _n^A:\mathbf{K}^a_n(Rep(A;B)) \rightarrow
\mathbf{K}^a_n(Rep(A\otimes \mathcal{K};B))$ be the homomorphism induced by
the functor $\varepsilon^A$. It is easy to verify that the family
$\{\varepsilon_n^A\}$ is a natural transformation from the functor
$\mathbf{K}^a_n(Rep(-;B))$ to $\mathbf{K}^a_n(Rep(-\otimes \mathcal{K};B))$.
Therefore, the following diagram
$$
\begin{CD}
  \mathbf{K}^a_n(Rep(A\otimes
\mathcal{K};B)) @ >e_n^{A}>>   \mathbf{K}^a_n(Rep(A;B))\\
  @V\varepsilon _n^{A\otimes \mathcal{K}}VV   @VV \varepsilon _n^AV\\
  \mathbf{K}^a_n(Rep(A\otimes
\mathcal{K}\otimes \mathcal{K};B))  @>>e_n^{A\otimes \mathcal{K}}>
\mathbf{K}^a_n(Rep(A\otimes
\mathcal{K};B))  \\
\end{CD}
$$
commutes, and it shows that for our purposes it suffices to verify that
\begin{equation}
\label{retr} e^{A}_n\varepsilon ^A_n=\id_{\mathbf{K}^a_n(Rep(A;B))},
\end{equation}
for all $A\in S_G$. Indeed, first note that

1. The equality \ref{retr} shows that $e^{A}_n$ is an epimorphism.

2. Since $\mathbf{K}^a_n(Rep(-\otimes\mathcal{K};B))$ is a stable functor,
according to the identity \ref{retr}, one easily shows that
$\varepsilon_n^{A\otimes\mathcal{K}}$ and $e_n^{A\otimes\mathcal{K}}$ are
isomorphisms (cf. Proposition 10.6 in \cite{suw}). Therefore $e^{A}_n$ is a
monomorphism.

\textit{Checking the equality} \ref{retr}. We construct a useful
isometry $\sigma_E:E\rightarrow E\otimes_k\mathcal{H}$, for any
countably generated $B$-module $E$, which is a morphism from
$(E,\phi)$ to $(E\otimes_k\mathcal{H},
(\phi\otimes\id_\mathcal{K})e_{A})$.

Choose $y\in\mathcal{H}$ so that $p(y)=y$ and $||y||=1$ and consider a
$B$-homomorphism $\sigma_E$ given by $x\mapsto x\otimes y$. For any $z\in
\mathcal{H}$ there exists $\lambda_z\in k$ determined uniquely by the
equation $p(z)=\lambda_zy$. Define $\sigma_E^*$ by the map $x\otimes z\mapsto
\lambda_zx$. The $B$-homomorphism $\sigma_E^*$ is adjoint to $\sigma_E$.
Since $\sigma_E^*\sigma_E(x)=\sigma_E^*(x\otimes y)=x$, one concludes that
$\sigma_E$ is an isometry. Since
$\sigma_E\phi(a)=((\phi\otimes\id_\mathcal{K})e_A(a))\sigma_E,$ the isometry
$\sigma_E$ is a morphism from $(E,\phi)$ into
$(E\otimes_k\mathcal{H},(\phi\otimes\id_\mathcal{K})e_A)$.

Consider restriction of $e^*_A\varepsilon$ to $D_\phi(A;E;B)$. We have a
$*$-homomorphism
\begin{equation}
(e^*_A\varepsilon )_E:D_\phi(A;E;B)\rightarrow
D_{(\phi\otimes\id)e_A}(A;E\otimes_k\mathcal{H};B)
\end{equation}
mapping $x\in D_\phi(A;E;B)$ to $x\otimes\id_{\mathcal{H}}\in
D_{(\phi\otimes\id)e_A}(A;E\otimes_k\mathcal{H};B)$. But
$$
(\sigma_Ex)(z)=(\sigma_E)(x(z))=x(z)\otimes y=((x\otimes
p)\sigma_E)(z)\textrm{ for any }x\in D_\phi(A;E;B), z\in E.
$$
Since $\sigma_E\sigma_E^*=\id_E\otimes p$, one concludes that
$\sigma_Ex\sigma^*_E=x\otimes p$. Therefore
\begin{equation}
\label{innerel} (e^*_A\varepsilon)_E(x)=\sigma_Ex\sigma_E^*+x\otimes(1-p).
\end{equation}
Let $\psi(x)=(e^*_A\varepsilon)_E(x)$, $\psi_0(x)=x\otimes p$ and
$\psi_1(x)=x\otimes(1-p)$. Then $\psi_0$  and $\psi_1$ are $*$-homomorphisms,
$\psi=\psi_0+\psi_1$ and
\begin{equation}
\label{ssss} \psi_0(x)\psi_1(x)=\psi_1(x)\psi_0(x)=0.
\end{equation}
Here we show that $e^*_A\varepsilon$ induces the identity homomorphism of the
group $\mathbf{K}_n^a(Rep(A,B))$ onto itself. Indeed, choose an element
$r\in\mathbf{K}_n^a(Rep(A,B))$. By definition of $\mathbf{K}^a_n$-groups the
element $r$ is represented by an element $r_\phi\in K_n^a(D_\phi(A;E;B))$.
Then the element $\mathbf{K}^a_n(e^*_A\varepsilon)(r)$ is represented by the
element
\begin{equation}
\label{rett}
K^a_n((e^*_A\varepsilon)_E)(r_\phi)=K^a_n(\psi_0+\psi_1)(r_\phi).
\end{equation}
Since $K^a_n$ is an additive functor, according to \ref{ssss} and
Lemma 2.1.18 in (\cite{hig3}), it follows that
$$
K^a_n(\psi_0+\psi_1)(r_\phi)=K^a_n(\psi_0)(r_\phi)+K^a_n(\psi_1)(r_\phi).
$$
Since $\psi_0=ad(s_E)$, the class of $K^a_n(\psi_0)(r_\phi)$ is equal to the
class of $r_\phi$. Thus the proof will be completed if we show that the class
of $K^a_n(\psi_1)(r_\phi)$ in $\mathbf{K}^a_n(Rep(A;B))$ is zero. Indeed, let
$$
s: E\otimes\mathcal{H}\rightarrow(E\otimes\mathcal{H})\oplus\mathcal{H}_B)
$$
be the isometry defined by $e\otimes h\mapsto e\otimes h\oplus 0$. Then the
$*$-homomorphism $ad(s)\psi_1$ may be factored through the $*$-homomorphism
$$
D_\phi(A;E;B)\rightarrow D_0(A;E\otimes(1-p)\mathcal{H}\oplus
\mathcal{H}_B;B).
$$
But, according to Kasparov's stabilization theorem, one has
$$
D_0(A;E\otimes (1-p)\mathcal{H}\oplus \mathcal{H}_B; B)\approx
\mathcal{L}(\mathcal{H}_B).
$$
Since $K^a_n(\mathcal{L}(\mathcal{H}_B))=0$, one concludes that classes of
$K^a_n(\psi_1)(z)$ and $K^a_n(ad(s)(\psi_1))(z)$ are equal to zero in
$\mathbf{K}^a_n(Rep(A;B))$. Thus the homomorphism
$\mathbf{K}^a_n(e_A\varepsilon)$ is the identity.
\end{proof}

\section{On the Isomorphism $\mathbf{K}^a_0(\mathrm{Rep}(-;B))\simeq KK_{-1}(-;B)$}

First we recall the definition of Kasparov's group $K^1(A,B)$, which will be
denoted by $KK_{-1}(A,B)$, where $A$ and $B$ are trivially graded
$C^*$-algebras with actions of a second countable compact group.

Consider a triple $(\varphi,E;p)$, where $E$ is a trivially graded countably
generated right $B$-module, $\varphi:A\rightarrow\mathcal{L}_B(E)$ is a
$*$-homomorphism and $p\in\mathcal{L}_B(E)$ is an invariant element so that
\begin{equation}
\label{kas}
\begin{array}{c}
p\varphi (a)-\varphi (a)p\in \mathcal{K}_B(E),\\
(p^{*}-p)\varphi(a)\in\mathcal{K}_B(E),\;\;(p^2-p)\varphi(a)\in
\mathcal{K}_B(E)
\end{array}
\end{equation}
for all $a\in A$. Such a triple will be called a \emph{Kasparov-Fredholm
$A,B$- module}. If all left parts in \ref{kas} are zero, then such a triple
is said to be \textit{degenerate}.

Define the sum of Kasparov-Fredholm $A,B$-modules by the formula
$$
(\varphi,E;p)\oplus(\varphi',E';p')=(\varphi\oplus\varphi',E\oplus E';p\oplus
p').
$$
Consider the equivalence relations:
\begin{itemize}
\item({\em Unitary isomorphism})
$A,B$-modules $(\varphi,E;p)$ and $(\varphi',E';p')$ will be said to be
unitarily isomorphic if there exists a unitary isomorphism $u:E\rightarrow
E'$ such that
$$
u\varphi(a)u^*=\varphi'(a),\;\;upu^*=p'
$$
for all $a\in A$.
\item({\em Homology})
$A,B$-modules $(\varphi,E;p)$ and $(\varphi',E;p')$ will be said to be
homologous if
$$
p'\varphi'(a)-p\varphi(a)\in\mathcal{K}_B(E)
$$
for all $a\in A$.
\end{itemize}
Simple checking shows that the equivalence relations defined above are well
behaved with respect to sum.

Let $\mathcal{E}^1(A,B)$ be the abelian monoid of classes of $A,B$-modules
with respect to the equivalence relation generated by the unitary isomorphism
and homology. Denote by $\mathcal{D}^1(A,B)$ the submonoid of
$\mathcal{E}^1(A,B)$ consisting of only those classes which are classes of
all degenerate triples. By definition
$$
E^1(A,B)=\mathcal{E}^1(A,B)/\mathcal{D}^1(A,B).
$$
Using the Kasparov stabilization theorem, one easily shows that the
definition of $E^1(A,B)$ coincides with Kasparov's original definition of
$E^1(A,B)$ which is isomorphic to $KK_{-1}(A,B)$ by lemma 2 of section 7 of
\cite{kas2}.

Recall that objects of $\mathrm{Rep}(A,B)$, by definition, have form
$(\varphi,E;p)$, where $p:(\varphi,E)\rightarrow(\varphi,E)$ is a projection
in the category $Rep(A;B)$. More precisely,
$$
\varphi(a)p-p\varphi(a)\in\mathcal{K}_B(E),\;\;p^*=p,\;\;p^2=p.
$$
A unitary isomorphism $s:(\varphi,E;p)\rightarrow(\psi,E',q)$ in
$\mathrm{Rep}(A,B)$ is a usual partial isometry $s:E\rightarrow E'$ such that
$$
s\varphi(a)-\psi(a)s\in\mathcal{K}_B(E,E'),\;\;s^*s=p,\;\;ss^*=q.
$$
Let $\tilde{\mathcal{E}}^1(A,B)$ be the abelian monoid of unitary isomorphism
classes of objects in $\mathrm{Rep(A;B)}$. According to Lemma \ref{lemaux},
one easily checks that the Grothendieck group of $\tilde{\mathcal{E}}^1(A,B)$
may be identified with $K_0(\mathrm{Rep}(A,B))$ (cf. \cite{kan1}).

There is a natural homomorphism
$$
\lambda _1:K_0(\mathrm{Rep}(A,B))\rightarrow E^1(A,B),
$$
defined by the map $(\varphi,E;p)\mapsto(\varphi,E;p)$. Indeed, let
$s:(\varphi,E;p)\rightarrow(\varphi',E';p')$ be a unitary isomorphism in
$\mathrm{Rep(A;B)}$. Consider the isomorphism in $Rep(A;B)$
$$
\bar s:(\varphi\oplus\psi, E\oplus E')\rightarrow(\psi\oplus\varphi,E'\oplus
E),
$$
where
$$
\bar s=\left(
\begin{array}{cc}
s&1-ss^*\\
1-s^*s&s
\end{array}
\right).
$$
It is clear that $(\varphi\oplus\psi,E\oplus E',\bar p)$ is isomorphic to
$(s(\varphi\oplus\psi)s^*,E'\oplus E,\bar q)$, which is homologous to
$(\psi\oplus\varphi,E'\oplus E,\bar q)$, with
$$
\bar p=\left(
\begin{array}{cc}
p&0\\
0&1
\end{array}
\right)\;\;\;\textrm{and}\;\;\; \bar p'=\left(
\begin{array}{cc}
1&0\\
0&p'
\end{array}
\right).
$$
This means that classes of $(\varphi,E;p)$ and $(\varphi',E';p')$ coincide in
$E^1(A,B)$.

Let $\mathrm{Rep}_J(A,B)=\mathrm{Rep}(A,B)/\mathrm{D}(A,J;B)$ be the
idempotent-complete $C^*$-category universally obtained from the category
$Rep(A,B)/D(A,J;B))$ (see section \ref{exe}). Let $(\varphi;E;p)$ be a
Kasparov-Fredholm $A,B$-module. Then $p$ defines a projector $\dot p$ in the
category $Rep_A(A,B)$. Thus the triple $(\varphi;E;\dot p)$ is an object in
$\mathrm{Rep}_A(A,B)$.

There is a well-defined homomorphism
$$
\mu:E^1(A,B)\rightarrow K_0(\mathrm{Rep}_A(A,B))
$$
defined by $(\varphi;E;p)\mapsto(\varphi;E;\dot p)$. This is checked below.

We recall definition of operatorial homotopy:
\begin{itemize}
\item({\em Operatorial homotopy}) An $A,B$-module $(\varphi,E;p)$ is
operatorially homotopic to $(\varphi,E;p')$ if there exists a continuous map
$p_t:[0;1]\rightarrow\mathcal{L}_B(E)$ such that $(\varphi,E;p_t)$ is an
$A,B$-module for any $t\in[0;1]$.
\end{itemize}

If $(\varphi,E;p)$ is homologous to $(\psi,E;q)$, then
$(\varphi,E;p)\oplus(\psi,E;0)$ is operatorially homotopic to
$(\varphi,E;0)\oplus(\psi,E;q)$. Indeed, the desired homotopy is defined by
the formula
$$
\left(\left(
\begin{array}{cc}
\varphi&0\\
0&\psi
\end{array}
\right),E\oplus E,\frac1{1+t^2}\left(
\begin{array}{cc}
p&tpq\\
tqp&t^2q
\end{array}
\right)\right),\;\;\;t\in[0;\infty]
$$
(cf. section 7 in \cite{kas2}). Thus the projections $\dot{p\oplus 0}$ and
$\dot{0\oplus q}$ are homotopic. Then, using Lemma 4 from section 6 in
\cite{kas2}, one concludes that the objects $(\varphi,E;\dot
p)\oplus(\psi,E;\dot0)$ and $(\varphi,E;\dot0)\oplus(\psi,E;\dot q)$ are
unitarily isomorphic objects in $\mathrm{Rep}_A(A,B)$. Let $(\varphi,E;p)$ be
unitarily isomorphic to $(\psi,E;q)$. Then $(\varphi,E;\dot p)$ is isomorphic
to $(\psi,E;\dot q)$ in the category $\mathrm{Rep}_A(A,B)$. Therefore $\mu$
is well-defined.

We are now ready to prove the following theorem.

\begin{theorem}
\label{iginol}
The natural homomorphism
$$
\lambda_1:K_0(\mathrm{Rep}(A,B))\rightarrow E^1(A,B)\simeq KK_{-1}(A,B)
$$
is an isomorphism.
\end{theorem}

\begin{proof} The homomorphism $\lambda_1$ is an epimorphism. Indeed,
let $(\varphi,E;p)$ be a Kasparov-Fredholm $A,B$-module. Applying techniques
of the Lemmata 17.4.2-17.4.3 in \cite{bla}, one can suppose that $p^*=p$ and
$||p||\leq1$. Then it is equivalent to $(\varphi\oplus0,E\oplus E;p')$, where
$$
p'=\begin{pmatrix}
p&\sqrt{p-p^2}\\
\sqrt{p-p^2}&1-p
\end{pmatrix}.
$$
Simple checking shows that $p'$ is a projection and $(\varphi\oplus0,E\oplus
E;p')$ is an object in $\mathrm{Rep}(A,B)$. To show that $\lambda_1$ is a
monomorphism, consider the commutative diagram
$$
\begin{CD}
K_0(\mathrm{Rep}(A,B))@>\lambda_1>>E^1(A,B)\\
@V\|VV                             @VV\mu V \\
K_0(\mathrm{Rep}(A,B))@>\xi>>K_0(\mathrm{Rep}_A(A,B)).
\end{CD}
$$
By Theorem \ref{topic}, $\xi$ is an isomorphism. Therefore $\lambda_1$ is a
monomorphism.
\end{proof}

\begin{center}
{\bf Acknowledgements}
\end{center}

\;

The work has synthesized ideas of many mathematicians. The author
would like to thank, especially those of, Cuntz, Higson, Karoubi,
Kasparov, Pedersen, Quillen, Rosenberg, Suslin, Weibel, Wodzicki.

It is great pleasure to acknowledge to J. Rosenberg and the
referee for the valuable comments.

The author was partially supported by INTAS grant 00 566, INTAS
grant 03-51-3251 and GRDF grant GEM1-3330-TB-03.


\end{document}